\pdfoutput=1
\documentclass{amsart}

\usepackage{amssymb}
\usepackage[english,francais]{babel}

\newtheorem{theorem}{Theorem}[section]

\newtheorem{e-proposition}[theorem]{Proposition}
\newtheorem{corollary}[theorem]{Corollary}
\newtheorem{e-definition}[theorem]{Definition}
\newtheorem{remark}[theorem]{Remark}

\newtheorem{theoreme}{Th\'eor\`eme}[section]

\newtheorem{proposition}[theoreme]{Proposition}

\def\og{\leavevmode\raise.3ex\hbox{$\scriptscriptstyle\langle\!\langle$~}}
\def\fg{\leavevmode\raise.3ex\hbox{~$\!\scriptscriptstyle\,\rangle\!\rangle$}}

\usepackage{amsmath}
\usepackage{mathrsfs}
\usepackage[all]{xy}
\usepackage{color}
\definecolor{red}{rgb}{1,0,0}

\definecolor{orange}{RGB}{255,127,0}

\DeclareMathOperator{\Hom}{Hom}
\DeclareMathOperator{\End}{End}
\DeclareMathOperator{\id}{id}
\DeclareMathOperator{\PBW}{PBW}

\DeclareMathOperator{\LC}{LC}
\DeclareMathOperator{\hol}{hol}
\newcommand{\abs}[1]{\left| #1 \right|} % absolute value
\newcommand{\crochet}[2]{\left[ #1 , #2 \right]} % Lie bracket
\newcommand{\into}{\hookrightarrow}

\newcommand{\NN}{\mathbb{N}} 
\newcommand{\ZZ}{\mathbb{Z}} 
\newcommand{\RR}{\mathbb{R}} 
\newcommand{\CC}{\mathbb{C}} 
\newcommand{\Koszulsign}[1]{\varepsilon(#1)}
\newcommand{\cinf}[1]{C^{\infty}(#1)}
\newcommand{\shuffle}[2]{\mathfrak{S}_{#1}^{#2}}
\newcommand{\sections}[1]{\Gamma(#1)}
\newcommand{\enveloping}[1]{\mathcal{U}(#1)}
\newcommand{\graded}{^{\bullet}}
\newcommand{\zo}{^{0,1}}
\newcommand{\oz}{^{1,0}}

\newcommand{\OO}{{\Omega}}
\newcommand{\sheaf}[1]{\mathscr{#1}}
\newcommand{\groupoid}[1]{\mathscr{#1}}
\newcommand{\category}[1]{\mathcal{#1}}
\newcommand{\lie}[2]{[#1,#2]} % Lie bracket

\begin{document}
\selectlanguage{english}

\title[Exponential map and $L_\infty$ algebra associated to a Lie pair]
{Exponential map and $L_\infty$ algebra \\ associated to a Lie pair}
\thanks{Research partially supported by the National Science Foundation [DMS-1101827].}

\author{Camille Laurent-Gengoux}
\address{D\'epartement de math\'ematiques, universit\'e de Lorraine}
\email{camille.laurent-gengoux@univ-lorraine.fr}
\author{Mathieu Sti\'enon}
\address{Department of Mathematics, Penn State University}
\email{stienon@math.psu.edu} 
\author{Ping Xu}
\address{Department of Mathematics, Penn State University}
\email{ping@math.psu.edu}

\maketitle

\begin{center}
{D\'edi\'e \`a la m\'emoire de Paulette Libermann.}
\end{center}

\begin{abstract}
\selectlanguage{english}% Text of abstract in English
In this note, we unveil homotopy-rich algebraic structures generated by the Atiyah classes relative to a 
Lie pair $(L,A)$ of algebroids. In particular, we prove that the quotient $L/A$ of such a pair admits an essentially 
canonical homotopy module structure over the Lie algebroid $A$, which we call Kapranov module. 

{\it Camille Laurent-Gengoux, Mathieu Sti\'enon, Ping Xu, C. R. Acad. Sci. Paris, Ser. I 350 (2012) 817-821.}

%\vskip 0.5\baselineskip
%
%\selectlanguage{francais}
%% Text of abstract in French
%\noindent{\bf R\'esum\'e} \vskip 0.5\baselineskip \noindent
%{\bf Application exponentielle et alg\`ebre $L_\infty$ associ\'ee \`a une paire de Lie}
%Dans cette note, nous d\'evoilons des structures alg\'ebriques, riches en homotopies, 
%engendr\'ees par les classes d'Atiyah relatives \`a une paire de Lie $(L,A)$ d'alg\'ebro\"ides. 
%En particulier, nous prouvons que le quotient $L/A$ d'une telle paire admet 
%une structure essentiellement canonique de module \`a homotopie pr\`es sur l'alg\'ebro\"ide de Lie $A$ 
%que nous appelons module de Kapranov. 
%{\it Pour citer cet article~: A. Name1, A. Name2, C. R. Acad. Sci.
%Paris, Ser. I 340 (2005).}

\end{abstract}

\selectlanguage{english}

\section{Kapranov modules}\label{AtiyahKapranov}

Let $A$ be a Lie algebroid (either real or complex) over a manifold $M$ with anchor $\rho$.
By an $A$-module, we mean a module of the corresponding Lie-Rinehart algebra $\sections{A}$ 
over the associative algebra $\cinf{M}$. 
An $A$-connection on a smooth vector bundle $E$ over $M$ is a bilinear map 
$\nabla:\sections{A}\otimes\sections{E}\to\sections{E}$ satisfying 
$\nabla_{fa} e=f\nabla_a e$ and
 $\nabla_a (fe)=\big(\rho(a)f\big) e+f\nabla_a e$, 
for all $a\in\sections{A}$, $e\in\sections{E}$, and $f\in\cinf{M}$. 
A vector bundle $E$ endowed with a \emph{flat} $A$-connection (also known as an infinitesimal $A$-action) 
is an $A$-module; more precisely, its space of smooth sections $\sections{E}$ is one.

\paragraph*{\textbf{Atiyah class}}
Given a Lie pair $(L,A)$ of algebroids, \emph{i.e.} a Lie algebroid $L$ with a Lie subalgebroid $A$, 
the Atiyah class $\alpha_E$ of an $A$-module $E$ relative to the pair $(L,A)$ is defined as the obstruction to the existence of 
an \emph{$A$-compatible} $L$-connection on $E$. An $L$-connection $\nabla$ is $A$-compatible if 
its restriction to $\sections{A}\otimes\sections{E}$ is the given infinitesimal $A$-action on $E$ and 
$\nabla_a\nabla_l-\nabla_l\nabla_a=\nabla_{[a,l]}$ for all $a\in\sections{A}$ and $l\in\sections{L}$. 
This fairly recently defined class (see~\cite{CSX}) 
has as double origin, which it generalizes, the Atiyah class of holomorphic vector bundles and the Molino class of foliations. 
The quotient $L/A$ of the Lie pair $(L,A)$ is an $A$-module~\cite{CSX}. Its Atiyah class $\alpha_{L/A}$ 
can be described as follows. Choose an $L$-connection $\nabla$ on $L/A$ extending the $A$-action. 
Its curvature is the vector bundle map $R^\nabla:\wedge^2 L\to\End(E)$
defined by $R^\nabla(l_1,l_2)=\nabla_{l_1}\nabla_{l_2}-\nabla_{l_2}\nabla_{l_1}
-\nabla_{\lie{l_1}{l_2}}$, for all $l_1, l_2\in\sections{L}$.
Since $L/A$ is an $A$-module, $R^\nabla$ vanishes on $\wedge^2 A$ and, therefore, 
determines a section $R^\nabla_{L/A}$ of $A^*\otimes(L/A)^*\otimes\End(L/A)$. 
It was proved in~\cite{CSX} that $R^\nabla_{L/A}$ is a $1$-cocycle for the Lie algebroid $A$ 
with values in the $A$-module $(L/A)^*\otimes\End(L/A)$ and that its
cohomology class $\alpha_{L/A}\in H^1\big(A;(L/A)^*\otimes\End(L/A)\big)$ 
is independent of the choice of the connection. 

\paragraph*{\textbf{Kapranov modules over a Lie algebroid}}
Let $M$ be a smooth manifold, and let $R$ be the algebra of smooth functions on $M$ valued in $\RR$ (or $\CC$). 
Let $A$ be a Lie algebroid over $M$. 
The Chevalley-Eilenberg differential $d_A$ and the exterior product make $\sections{\wedge\graded A^*}$ 
into a differential graded commutative $R$-algebra.
 
Now let $E$ be a smooth vector bundle over $M$. 
Deconcatenation defines an $R$-coalgebra structure on $\sections{S\graded E}$. 
%--->
The comultiplication 
$\Delta:\sections{S\graded E}\to\sections{S\graded E}\otimes_R\sections{S\graded E}$ is given by 
\[ \Delta(e_1\odot e_2\odot\cdots\odot e_n)=
\sum_{p+q=n}\sum_{\sigma\in\shuffle{p}{q}} (e_{\sigma(1)}\odot\cdots\odot e_{\sigma(p)}) \otimes
(e_{\sigma(p+1)}\odot\cdots\odot e_{\sigma(n)}) ,\] 
for any $e_1,\dots,e_n\in\sections{E}$. 
%<---
Let $\mathfrak{e}$ denote the ideal of $\sections{S\graded (E^*)}$ generated by $\sections{E^*}$. 
The algebra $\Hom_R\big(\sections{S\graded E},R\big)$ 
dual to the coalgebra $\sections{S\graded E}$ 
is the $\mathfrak{e}$-adic completion of $\sections{S\graded (E^*)}$. 
It will be denoted by $\sections{\hat{S}\graded (E^*)}$. 
Equivalently, one can think of the completion $\hat{S}\graded (E^*)$ of $S\graded (E^*)$ 
as a bundle of algebras over $M$. 
Note that $\sections{\wedge\graded A^*}$ is an $R$-subalgebra of 
$\sections{\wedge\graded A^*\otimes \hat{S}\graded E^*}$. 

Recall that an $L_\infty[1]$ algebra is a $\ZZ$-graded vector space 
$V=\bigoplus_{n\in\ZZ}V_n$ endowed with 
a sequence $(\lambda_k)_{k=1}^\infty$ of symmetric multilinear maps 
$\lambda_k: \otimes^k V\to V$ of degree $1$
satisfying the generalized Jacobi identity
\[ \textstyle\sum_{k=1}^n\sum_{\sigma\in\shuffle{k}{n-k}}\Koszulsign{\sigma;v_1,\cdots,v_n} \;
\lambda_{1+n-k}\big(\lambda_k(v_{\sigma(1)},\cdots,v_{\sigma(k)}),
v_{\sigma(k+1)},\cdots,v_{\sigma(n)}\big)=0 \]
for each $n\in\NN$ and for any homogeneous vectors $v_1,v_2,\dots,v_n\in V$. 
Here $\shuffle{p}{q}$ denotes the set of 
\emph{$(p,q)$-shuffles}\footnote{A \emph{$(p,q)$-shuffle} is a permutation $\sigma$ 
of the set $\{1,2,\cdots,p+q\}$ such that $\sigma(1)\le\sigma(2)\le\cdots\le\sigma(p)$ 
and $\sigma(p+1)\le\sigma(p+2)\le\cdots\le\sigma(p+q)$.}
and $\Koszulsign{\sigma; v_1, \cdots, v_n}$ the \emph{Koszul sign}\footnote{The \emph{Koszul sign} 
of a permutation $\sigma$ of the (homogeneous) vectors $v_1,v_2,\dots,v_n$ is determined by 
the relation $v_{\sigma(1)}\odot v_{\sigma(2)}\odot\cdots\odot v_{\sigma(n)} 
= \Koszulsign{\sigma; v_1, \cdots, v_n}\; v_1\odot v_2\odot\cdots\odot v_n$.} 
of the permutation $\sigma$ of the (homogeneous) vectors $v_1,v_2,\dots,v_n$.  

\begin{e-definition}\label{kapmod}
A Kapranov module over a Lie algebroid $A\to M$ is a vector bundle $E\to M$ 
together with an $L_\infty[1]$ algebra structure on $\sections{\wedge\graded A^*\otimes E}$ 
defined by a sequence $(\lambda_k)_{k\in\NN}$ of multibrackets (called Kapranov multibrackets) 
such that (1) the unary bracket 
$\lambda_1:\sections{\wedge\graded A^*\otimes E}\to\sections{\wedge^{\bullet+1} A^*\otimes E}$
is the Chevalley-Eilenberg differential associated to an infinitesimal $A$-action on $E$, 
and (2) all multibrackets 
$\lambda_k:\otimes^k\sections{\wedge^{\bullet} A^*\otimes E}\to\sections{\wedge^{\bullet} A^*\otimes E}[1]$
with $k\geq 2$ are $\sections{\wedge\graded A^*}$-multilinear.
\end{e-definition}

\begin{proposition}\label{pro:Kmodule}
Let $A$ be a Lie algebroid over a smooth manifold $M$ and let $E$ be a smooth vector bundle over $M$. 
Each of the following four data is equivalent to a Kapranov $A$-module structure on $E$. 
\begin{enumerate}
\item \label{c1} 
A degree 1 derivation $D$ of the graded algebra $\sections{\wedge\graded A^*\otimes\hat{S}(E^*)}$, 
%which maps $\sections{\wedge^k A^*\otimes\hat{S}(E^*)}$ to $\sections{\wedge^{k+1} A^*\otimes\hat{S}(E^*)}$, 
which preserves the filtration $\sections{\wedge A^*\otimes\hat{S}^{\geq n}(E^*)}$, satisfies $D^2=0$, 
and whose restriction to $\sections{\wedge\graded A^*}$ is the Chevalley-Eilenberg differential 
of the Lie algebroid $A$. 
(Here, by convention, all elements of $\hat{S}(E^*)$ have degree 0.)
\item 
An infinitesimal action of $A$ 
on $\hat{S}(E^*)$ by derivations which preserve the decreasing filtration $\hat{S}^{\geq n}(E^*)$. 
\item 
An infinitesimal action of $A$ on $S(E)$ 
by coderivations which preserve \linebreak $S^{\geq 1}(E)$ and the increasing filtration $S^{\leq n}(E)$.
\item \label{c4} 
An infinitesimal action of $A$ on $E$ together with a sequence of morphisms of vector bundles 
$\boldsymbol{R}_k:S^k(E)\to A^*\otimes E$ ($k\geq 2$) whose sum 
\[ \boldsymbol{R}=\sum_{k=2}^\infty\boldsymbol{R}_k\in\sections{A^*\otimes\hat{S}(E^*)\otimes E} \] 
is a solution of the Maurer-Cartan equation $d_A\boldsymbol{R}+\tfrac{1}{2}\crochet{\boldsymbol{R}}{\boldsymbol{R}}=0$. 
(Here, we consider $\sections{\hat{S}(E^*)\otimes E}$ as the space of formal vertical vector fields on $E$ 
along the zero section and derive a natural Lie bracket on the graded vector space 
$\sections{\wedge^{\bullet}A^*\otimes\hat{S}(E^*)\otimes E}$.)
\end{enumerate}
\end{proposition}

Characterizations (\ref{c1}) and (\ref{c4}) are related by the identity $D=d_A^{\hat{S}(E^*)}+\boldsymbol{R}$, 
where $d_A^{\hat{S}(E^*)}$ denotes the Chevalley-Eilenberg differential associated to the infinitesimal $A$-action on $E$, 
and $\boldsymbol{R}$ denotes its own action on $\sections{\wedge\graded A^*\otimes \hat{S}(E^*)}$ by contraction.
On the other hand, for any $k\geq 2$, the $k$-th Kapranov multibracket $\lambda_k$ is related to the $k$-th component 
$\boldsymbol{R}_k\in\sections{A^*\otimes S^k E^* \otimes E}$ of the Maurer-Cartan element $\boldsymbol{R}$ 
through the equation 
\[ \lambda_k(\xi_1\otimes e_1,\cdots,\xi_k\otimes e_k)=(-1)^{\abs{\xi_1}+\cdots+\abs{\xi_k}}
\xi_1\wedge\cdots\wedge\xi_k\wedge\boldsymbol{R}_k (e_1,\cdots,e_k) ,\]
which is valid for any $e_1,\dots,e_k\in\sections{E}$ 
and any homogeneous elements $\xi_1,\dots,\xi_k$ of $\sections{\wedge\graded A^*}$.

The algebraic structure described in the above proposition is related to 
Costello's $L_\infty$ algebras over the differential graded algebra 
$(\sections{\wedge\graded A^*},d_A)$~\cite{arXiv:1112.0816}, 
and to Yu's $L_\infty$ algebroids~\cite{Yu}.

Two Kapranov $A$-modules $E_1$ and $E_2$ over $M$ are isomorphic
 if there exists an isomorphism $\Phi:S(E_1)\to S(E_2)$
of bundles of coalgebras over $M$, which intertwines 
the infinitesimal $A$-actions. 

\section{Exponential map and Poincar\'e-Birkhoff-Witt isomorphism}

Assume $\groupoid{A}$ is a Lie subgroupoid of a Lie groupoid $\groupoid{L}$
(over the same unit space), 
and let $A$ and $L$ denote the corresponding Lie algebroids.  
The source map $s:\groupoid{L}\to M$ factors through the quotient 
of the action of $\groupoid{A}$ on $\groupoid{L}$ by multiplication from the right. 
Therefore, it induces a surjective submersion $s:\groupoid{L}/\groupoid{A}\to M$. 
Note that the zero section $0:M\to L/A$ and the unit section $1:M\to \groupoid{L}/\groupoid{A}$ are both embeddings of $M$. 

\begin{proposition}\label{expo}
Each choice of a splitting of the short exact sequence of vector bundles $0\to A\to L\to L/A\to 0$ 
and of an $L$-connection $\nabla$ on $L/A$ extending the $A$-action 
determines an exponential map, \emph{i.e.} a fiber bundle map 
\[ \exp^\nabla:L/A\to\groupoid{L}/\groupoid{A} ,\]
%\[ \xymatrix{ L/A \ar[r]^{\exp^\nabla} \ar[d]_{\pi} & \groupoid{L}/\groupoid{A} \ar[d]^s \\  M \ar[r]_{\id} & M , } \]
which identifies the zero section of $L/A$ to the unit section of $\groupoid{L}/\groupoid{A}$, 
whose differential along the zero section of $L/A$ is the canonical isomorphism 
between $L/A$ and the tangent bundle to the $s$-foliation of $\groupoid{L}/\groupoid{A}$ along the unit section, 
and which is locally diffeomorphic around $M$. 
\end{proposition}

Let $\mathcal{N}(L/A)$ denote the space of all functions on $L/A$ which, together with their derivatives 
of all degrees in the direction of the $\pi$-fibers, vanish along the zero section. 
The space of $\pi$-fiberwise differential operators on $L/A$ along the zero section 
is canonically identified to the symmetric $R$-algebra $\sections{S(L/A)}$.
Therefore, we have the short exact sequence of $R$-algebras 
\begin{equation}\label{eq:jet1}
 0 \to \mathcal{N}(L/A) \to \cinf{L/A} \to \Hom_R\big(\sections{S(L/A)},R\big) \to 0 .
\end{equation}

Likewise, let $\mathcal{N}(\groupoid{L}/\groupoid{A})$ denote the space of all functions 
on $\groupoid{L}/\groupoid{A}$ which, together with their derivatives of all degrees in the direction of the $s$-fibers, 
vanish along the unit section. 
The space of $s$-fiberwise differential operators on $\groupoid{L}/\groupoid{A}$ 
along the unit section is canonically identified to 
the quotient of the enveloping algebra $\mathcal{U}(L)$ by the left ideal generated by $\sections{A}$.
Therefore, we have the short exact sequence of $R$-algebras
\begin{equation}\label{eq:jet2}
 0 \to \mathcal{N}(\groupoid{L}/\groupoid{A}) \to \cinf{\groupoid{L}/\groupoid{A}} \to 
\Hom_R\big(\tfrac{\mathcal{U}(L)}{\mathcal{U}(L)\sections{A}},R\big) \to 0 .
\end{equation}

Since the exponential (or more precisely its dual)  
maps $\mathcal{N}(\groupoid{L}/\groupoid{A})$ to $\mathcal{N}(L/A)$, 
it induces an isomorphism of $R$-modules from $\Hom_R\big(\sections{S(L/A)},R\big)$ to \linebreak 
$\Hom_R\big(\tfrac{\mathcal{U}(L)}{\mathcal{U}(L)\sections{A}},R\big)$.

\begin{proposition}\label{exp2pbw}
Each choice of a splitting of the short exact sequence of vector bundles
$0\to A\to L\to L/A\to 0$ and of an $L$-connection $\nabla$ on $L/A$ extending the $A$-action 
determines an isomorphism of filtered $R$-modules 
\[ \PBW:\sections{S(L/A)}\to\tfrac{\mathcal{U}(L)}{\mathcal{U}(L)\sections{A}} \] 
called Poincar\'e-Birkhoff-Witt map. 
\end{proposition}

\begin{remark}
In case $L=A\bowtie B$ is the Lie algebroid sum of a matched pair of Lie algebroids $(A, B)$, 
the $L$-connection $\nabla$ on $L/A\cong B$ extending the $A$-action determines a $B$-connection on $B$, 
the coalgebras $\tfrac{\enveloping{L}}{\enveloping{L}\sections{A}}$ and $\enveloping{B}$ are isomorphic, 
and the corresponding Poincar\'e-Birkhoff-Witt map $\PBW:\sections{S(B)}\to\enveloping{B}$ 
is standard (see~\cite{NWX} for instance).
\end{remark}

\begin{proposition}\label{recursive}
The Poincar\'e-Birkhoff-Witt map associated to a splitting $j:L/A \to L$ of the short exact sequence of vector bundles 
$0\to A\to L\to L/A\to 0$ and an $L$-connection $\nabla$ on $L/A$ 
satisfies $\PBW(1)=1$ and, for all $b\in\sections{L/A}$ and $n\in\NN$, $\PBW(b)=j(b)$ and 
$\PBW(b^{n+1})=j(b)\cdot\PBW(b^n)-\PBW\big(\nabla_{j(b)} (b^n)\big)$, 
where $b^k$ stands for the symmetric product $b\odot b\odot \cdots \odot b$ of $k$ copies of $b$.
\end{proposition}

\begin{remark}
Although the construction of the Poincar\'e-Birkhoff-Witt map outlined above presupposes 
that $L$ and $A$ are integrable real Lie algebroids, $\PBW$ can be defined for any real (resp.\ complex) 
Lie pair provided one works with local (resp.\ formal) groupoids. 
\end{remark}

%--->
Recall that $L/A$ is an $A$-module. The representation of $A$ on $L/A$, 
an obvious generalization of the Bott connection, induces an infinitesimal action of $A$ 
on the coalgebra $\sections{S(L/A)}$ by coderivation. 
Besides, multiplication by elements of $\sections{A}$ from the left in the coalgebra $\mathcal{U}(L)$ 
induces an infinitesimal action of $A$ on the coalgebra $\tfrac{\mathcal{U}(L)}{\mathcal{U}(L)\sections{A}}$ by coderivation.
%<--- 

The infinitesimal actions of $A$ on $L/A$ and $\groupoid{L}/\groupoid{A}$  
induce infinitesimal actions of $A$ by derivations on the algebras of functions $C^\infty (L/A)$ 
and $C^\infty (\groupoid{L}/\groupoid{A})$ and, consequently, on the algebras of infinite jets 
$\Hom_R\big(\sections{S(L/A)},R\big)$ and $\Hom_R\big(\tfrac{\mathcal{U}(L)}{\mathcal{U}(L)A},R\big)$.

\begin{proposition}\label{pro:jet}
(1) The space $\Hom_R\big(\tfrac{\mathcal{U}(L)}{\mathcal{U}(L)\sections{A}},R\big)$ 
of infinite $s$-fiberwise jets along $M$ of functions on $\groupoid{L}/\groupoid{A}$ 
is an associative algebra on which the Lie algebroid $A$ acts infinitesimally by derivations. 
(2) The dual of the exponential map 
$\PBW^*:\Hom_R\big(\tfrac{\mathcal{U}(L)}{\mathcal{U}(L)\sections{A}},R\big)\to \Hom_R\big(\sections{S (L/A)},R\big)$ 
is an isomorphism of associative algebras, which may or may not intertwine the infinitesimal $A$-actions.
\end{proposition}

\section{$L_\infty[1]$ algebra associated to a Lie pair}

Our main result is the following

\begin{theorem}\label{thm:main}
If $(L,A)$ is a Lie pair, \emph{i.e.} a Lie algebroid $L$ together with a Lie subalgebroid $A$, 
then $L/A$ admits a Kapranov module structure, canonical up to isomorphism, 
over the Lie algebroid $A$, whose $\boldsymbol{R}_2 \in \Gamma (A^* \otimes S^2 (L/A)^* \otimes L/A)$ 
(see Proposition~\ref{pro:Kmodule}) is a 1-cocycle representative of the Atiyah class of $L/A$ 
relative to the pair $(L,A)$.
\newline 
Moreover, when $L=A\bowtie B$ is the Lie algebroid sum of a matched pair $(A, B)$ of Lie algebroids 
and there exists a torsion free flat $B$-connection $\nabla$ on $B$, the components of the Maurer-Cartan element 
$\boldsymbol{R}$ satisfy the recursive formula $\boldsymbol{R}_{k+1}=\partial^\nabla\boldsymbol{R}_k$, 
where $\partial^\nabla$ denotes the covariant differential associated to the connection.
\end{theorem}

\noindent\textbf{Sketch of proof}\;
Choose a splitting of the short exact sequence of vector bundles $0\to A\to L\to L/A\to 0$ 
and an $L$-connection $\nabla$ on $L/A$ extending the $A$-action.
Identify $\sections{\hat{S}(L/A)^*}$ to $\Hom_R\big(\tfrac{\enveloping{L}}{\enveloping{L}\sections{A}},R\big)$
via the $\PBW$ map and pull back the infinitesimal $A$-action of the latter to the former. 
According to Proposition~\ref{pro:Kmodule}, the resulting $A$-action on $\sections{\hat{S}(L/A)^*}$
 by derivations determines a Kapranov $A$-module structure on $L/A$. 
Making use of Proposition~\ref{recursive}, one can check directly that 
$\boldsymbol{R}_2$ is a 1-cocycle representative of the Atiyah class $\alpha_{L/A}$. \qed

As immediate consequences, we recover the following results of~\cite{CSX}.

\begin{corollary}\label{cor1}
Given a Lie algebroid pair $(L,A)$, let $\mathcal{U}(A)$ denote the universal enveloping algebra of the Lie algebroid $A$
and let $\category{A}$ denote the category of $\mathcal{U}(A)$-modules. 
The Atiyah class of the quotient $L/A$ makes $L/A[-1]$ 
into a Lie algebra object in the derived category $D^b(\category{A})$.
\end{corollary}

\begin{corollary}\label{cor2}
Let $(L,A)$ be a Lie pair and let $\sheaf{C}$ be a bundle (of finite or infinite rank)
of associative commutative algebras on which $A$ acts by derivations.
There exists an $L_\infty [1]$ algebra structure on 
$\sections{\wedge^\bullet A^*\otimes L/A\otimes \sheaf{C}}$, 
canonical up to $L_\infty$ isomorphism. 
Moreover, $H^{\bullet-1}(A;L/A\otimes\sheaf{C})$ is a graded Lie algebra 
whose Lie bracket only depends on the Atiyah class of $L/A$.
\end{corollary}

For $\sheaf{C}=\CC$, the Lie bracket on the cohomology $H^{\bullet-1} (A, L/A)$ happens to be trivial.

\section{An example due to Kapranov}

Let $X$ be a K\"ahler manifold with real analytic metric. 
Recall that the eigenbundles $T\zo_X$ and $T\oz_X$ of the complex structure $J:T_X\to T_X$ ($J^2=-\id$) 
form a matched pair of Lie algebroids~\cite{LSX}.
Fix a point $x\in X$. The exponential map $\exp^{\LC}_x:T_x X\to X$ defined using the geodesics of the Levi-Civita 
connection $\nabla^{\LC}$ originating from the point $x$ needs not be holomorphic. 

However, Calabi constructed a holomorphic exponential map $\exp^{\hol}_x:T_x X\to X$ 
as follows~\cite{Calabi} (see also~\cite{BCOV}).
First, extend the Levi-Civita connection $\CC$-linearly to a $T_X\otimes\CC$-connection $\nabla^\CC$ on $T_X\otimes\CC$. 
Since $X$ is K\"ahler, $\nabla^{\LC}J=0$ and $\nabla^\CC$ restricts to a $T_X\otimes\CC$-connection on $T\oz_X$. 
It is easy to check that the induced $T\zo_X$-connection on $T\oz_X$ is the canonical infinitesimal $T\zo_X$-action 
on $T\oz_X$ --- a section of $T\oz_X$ is $T\zo_X$-horizontal iff it is holomorphic --- while the induced $T\oz_X$-connection 
$\nabla\oz$ on $T\oz_X$ is flat and torsion free. 
Now let $X'$ denote the manifold $X$ and let $X''$ denote $X$ with the opposite complex structure $-J$. 
The image of the diagonal embedding $X\into X'\times X''$ is totally real so $X'\times X''$ can be seen 
as a complexification of $X$. The restriction of $T_{X'\times X''}$ (resp.\ its subbundle $T_{X'}\times X''$) 
along the diagonal $X$ is precisely the complexified tangent bundle $T_X\otimes\CC$ (resp.\ its subbundle $T\oz_X$). 
(See~\cite{Weinstein} for a discussion on integration of complex Lie algebroids.) 
The analytic continuation of the $T\oz_X$-connection $\nabla\oz$ on $T\oz_X$ in a neighborhood of the diagonal 
is a holomorphic $T_{X'}\times X''$-connection on the Lie algebroid $T_{X'}\times X''$, whose exponential map 
$\exp^{\hol}_x$ at a diagonal point $(x,x)$ takes $T_x X'\times\{x\}$ (which is $(T\oz_X)_x$ or $T_x X$) 
into $X'\times\{x\}$ (which is $X$). 

Consider the Lie pair $(L=T_{X'\times X''},A=X'\times T_{X''})$, the corresponding Lie groupoids 
$\groupoid{L}=(X'\times X'')\times(X'\times X'')$ and $\groupoid{A}=X'\times(X''\times X'')$, and the associated quotients 
$L/A=T_{X'}\times X''$ and $\groupoid{L}/\groupoid{A}=(X'\times X')\times X''$. 
Calabi's holomorphic exponential map $\exp^{\hol}$ is indeed the restriction along the diagonal of the exponential map 
$\exp^{\nabla\oz}:L/A\to\groupoid{L}/\groupoid{A}$ associated to the $T_{X'}\times X''$-connection $\nabla\oz$ 
on the Lie algebroid $T_{X'}\times X''$ as described in Proposition~\ref{expo}. 

Taking the infinite jet of $\exp^{\hol}$, we obtain, as in Proposition~\ref{exp2pbw}, a Poincar\'e-Birkhoff-Witt map 
$\PBW^{\hol}:\sections{S (T\oz_X)}\to\enveloping{T\oz_X}$.
Then, pulling back the infinitesimal $T\zo_X$-action on $\enveloping{T\oz_X}$ to an infinitesimal $T\zo_X$-action 
by coderivations on $\sections{S(T\oz_X)}$, we obtain, as in Theorem~\ref{thm:main}, a Kapranov $T\zo_X$-module 
structure on $T\oz_X$. In this context, the tensors $R_n\in\OO\zo\big(\Hom(S^n T\oz_X,T\oz_X)\big)$ are the curvature 
$R_2\in\OO^{1,1}\big(\End(T\oz_X)\big)$ and its higher covariant derivatives. 
Hence we recover the following result of Kapranov: 

\begin{theorem}[\cite{Kapranov}]
The Dolbeault complex $\OO^{0,\bullet}(T\oz_X)$ of a K\"ahler manifold is an $L_\infty[1]$ algebra. 
For $n\geq 2$, the $n$-th multibracket $\lambda_n$ 
%\[ \lambda_n:\OO^{0,j_1}(T\oz_X)\otimes\cdots\otimes\OO^{0,j_n}(T\oz_X)\to\OO^{0,j_1+\cdots+j_n+1}(T\oz_X) \] 
is the composition 
\[ \OO^{0,j_1}(T\oz_X)\otimes\cdots\otimes\OO^{0,j_n}(T\oz_X)\to\OO^{0,j_1+\cdots+j_n}(\otimes^n T\oz_X)
\to\OO^{0,j_1+\cdots+j_n+1}(T\oz_X) \] 
of the wedge product with the map 
induced by $R_n\in\OO\zo\big(\Hom(\otimes^n T\oz_X,T\oz_X)\big)$, while $\lambda_1$ is the 
Dolbeault operator $\overline{\partial}:\Omega^{0,j}(T\oz_X)\to\Omega^{0,j+1}(T\oz_X)$.
\end{theorem}

\section*{Acknowledgements}
Camille Laurent-Gengoux likes to express his gratitude to Penn State for its hospitality. 

\bibliographystyle{amsplain}
\bibliography{biblio}

\end{document}